\documentclass[12pt]{article}
\usepackage{amsmath,amssymb,amsfonts,amscd}
\usepackage{mitpress}

\makeatletter  
\renewcommand{\subsection}{\@startsection{subsection}{2}{0mm}{\baselineskip}{-1em}{\normalfont\normalsize\itshape}}  
\makeatother 
\newtheorem{theorem}{Theorem}

\newtheorem{corollary}[theorem]{Corollary}

\newtheorem{lemma}[theorem]{Lemma}

\newtheorem{proposition}[theorem]{Proposition}
\newtheorem{remark}[theorem]{Remark}

\newdimen\dummy
\dummy=\oddsidemargin
\begin{document}

\title{A non standard approach of spectral theory of self-adjoint operators
(generalized Gelfand eigenvectors)}
\author{Fatma Karray Meziou \\
%EndAName
Universit\'{e} Pierre \& Marie Curie.\\
meziou@math.jussieu.fr}
\maketitle

\begin{abstract}
Let $u=\int_{-\infty }^{+\infty }\lambda dE_{\lambda }$ be a self-adjoint
operator in a Hilbert space $H$. Our purpose is to provide a non-standard
description of the spectral family $(E_{\lambda })$ and the generalized
Gelfand eigenvectors.
\end{abstract}

\section{The non standard framework}

We consider here the nonstandard analysis framework (see \cite{Gold}, for a
brief overview) and apply it to set universes : Let $A$ be a set of $\in $%
-independent elements ( i.e. a non empty set such as there is not any chain $%
x_{0}\in x_{1}\in \ldots \in x_{n}$ whose ends belong to $A$ ) and let us
associate to $A$ its superstructure $U_{A}$, i.e. the union of the sequence $%
(A_{n})$\ defined by $A_{0}=A$\ and $A_{n+1}=A_{n}\cup P^{\prime }(A_{n})$,
where $P^{\prime }$\ denotes the set of non empty subsets. The couple made
of $U_{A}$\ and the restriction of $\in $\ to $U_{A}$\ is a standard model
of Zermelo theory with atoms and where Zorn's theorem holds.

\noindent \textbf{Note :}

\begin{enumerate}
\item Zermelo theory is ZF without the substitution axiom (see \cite{Jech}).

\item We need to choose once for all a $\varnothing _{A}\in A$ that will
play the role of the empty set in $U_{A}$; $U_{A}$ atoms are then the
elements of $A\backslash \{\varnothing _{A}\}.$ Such a $U_{A}$\ universe
contains all necessary sets to describe the structures (in the sense of
N.Bourbaki) over elements of $U_{A}$.

\item We may assume that $A$ is chosen such that $\mathbb{C}$ and $H$ are
its subsets; even if this requires to replace $\mathbb{C}$ and $H$ by
isomorphic copies. $\mathbb{C}$ and $H$ have then their binary operation
graphs in $U_{A}$. One can associate to $U_{A}$ a first order language $%
L(U_{A})$\ whose relational symbols are $=$ and $\in $ and whose constant
set is $U_{A}$.
\end{enumerate}

The non-standard method relies on the existence of a function called
enlargement defined from $U_{A}$ to the set class as $a\mapsto ^{\ast }a$
with the following properties:

\begin{description}
\item i. $^{\ast }a=a$ for each atom,

\item ii. if $R$ is a relation of $L(U_{A})$ without any free variable and
where all quantifiers are limited - i.e. in the form of $\exists x\in t$ or $%
\forall x\in t$ - and if $^{\ast }R$\ is the relation obtained by replacing
every outstanding constant $a$ by its image $^{\ast }a$, $R$ and $^{\ast }R$
relations are either both true or both false (transfer principle),

\item iii. if $F\in U_{A}$ is a filter on a set $u\in U_{A}$, the
intersection of\ $^{\ast }F$ sets (which is a filter over $^{\ast }u$) is
not empty.
\end{description}

\noindent In particular, given i and ii, $\mathbb{C}$ could be enlarged into
a set $^{\ast }\mathbb{C}$ containing $\mathbb{C}$. We will call $^{\ast }%
\mathbb{C}$ elements hypercomplexes and define over $^{\ast }\mathbb{C}$ an
addition and a multiplication operations, whose graphs are the enlargements
of the addition and multiplication graphs over $\mathbb{C}$. For these
operations $^{\ast }\mathbb{C}$ is an algebraically closed field. The same
analysis applies to $\mathbb{R}$ and leads to the subfield $^{\ast }\mathbb{R%
}$ of $^{\ast }\mathbb{C}$. $^{\ast }\mathbb{R}$ is totally ordered and
generates with element $i$ the field $^{\ast }\mathbb{C}$ (This allows to
consider the real and imaginary parts of a hypercomplex and to introduce the
absolute hypervalue of a hypercomplex).\newline

\noindent An element $\xi \in ^{\ast }\mathbb{R}$ is infinitesimal (resp.
limited), if $-1<n\xi<1$ for every natural number $n$ (resp. $-n<\xi <n$ for
at least one natural number). A hypercomplex is infinitesimal (resp.
limited) if its absolute hypervalue is infinitesimal (resp. limited). The
subset $^{\ast }\mathbb{C}_{b}$ made of bounded elements is a valuation ring
whose residual field could be identified with $\mathbb{C}$ : more precisely
the image of $z\in ^{\ast }\mathbb{C}_{b}$ is the unique complex number $%
^{\circ }z$ such that $z-^{\circ }z$ is infinitesimal.\newline

\noindent From now on, we consider a fixed universe $U$ and an enlargement $%
^{\ast }U$ of $U$.

\section{Hyperhermitian spaces}

\subsection{Internal vector spaces}

A set is called internal if it belongs to an enlargement. A $^{\ast }\mathbb{%
C}$-vector space $F$ is internal if the set $F$ is internal and if the
graphs of its binary operations are internal. A linear function $%
u:F\rightarrow G$ between internal vector spaces is internal if its graph is
internal.

\begin{theorem}
\label{thm1}Internal $^{\ast }\mathbb{C}$ vector spaces and internal linear
applications make up an abelian category (denoted $Ivs$).
\end{theorem}

Let $E$ be a complex vector space of $U$ and let us assume that binary
operations are in $U$. If $vss_{f}(E)$\ is the set of finite-dimensional
subspaces of $E$, its enlargement $^{\ast }vss(E)$\ is made of $^{\ast }E$
internal subspaces (and is itself an internal linear space for the
enlargement of binary operations of $E$). \newline
An internal vector space $F$ is said to be hyperfinite-dimensional, if there
is a complex vector space $E$ in $U$ such that $F$ is isomorphic (through an
internal isomorphism) to an element of $^{\ast }vss(E)$.

\begin{corollary}
\label{coro2} The full subcategory of $Ivs$ set made of hyperfinite-\newline
dimensional subspaces is a thick subcategory.
\end{corollary}

\noindent The importance of these subspaces lies in the following
proposition :

\begin{proposition}
\label{prop3} Let $H$ be a complex vector space of $U$. Since $\mathbb{C}%
\subset \, ^{\ast }\mathbb{C}$, $^{\ast }H$ is also a complex space and the
function : 
\begin{equation*}
\begin{tabular}{ll}
$f:$ & $H\rightarrow ^{\ast }H$ \\ 
& $h\mapsto ^{\ast }h$%
\end{tabular}%
\end{equation*}%
\ is an injective linear map whose image is included in an internal\newline
hyperfinite-dimensional subspace.
\end{proposition}

\subsection{Definition of hyperhermitian spaces and orthonormal hyperbases}

The concept of positive-definite sesquilinear form extends easily to
hypercomplex vector spaces. We call hyperhermitian space an internal
hyperfinite-dimensional vector space provided with an internal
positive-definite sesquilinear form. A number of properties are preserved by
transfer, when passing from hermitian to hyperhermitian spaces. Let us
mention in particular the orthogonal property and moreover the decomposition
of a hyperhermitian space into a direct sum of an internal linear subspace
and its orthogonal.

Orthonormal hyperbases play a very important role here : a subset $B$ of a
hyperhermitian space $F$ is an orthonormal hyperbasis, if $B$ is an
orthonormal set and $0$ is the only element of $F$ to be orthogonal to every
element of $B$.

\begin{theorem}
\label{thm4}Let $F$ be a hyperhermitian vector space with scalar product $%
<,> $ and let $u$ be a symmetric internal linear function for $<,>$ over $F$%
. There exists an orthonormal hyperbasis $B$ of $F$ made of eigenvectors of $%
u$.
\end{theorem}

Such a hyperbasis $B$ has an important property: it is hyperfinite i.e. it
belongs to the enlargement of the set of finite subsets of some set in $%
U_{A} $. We can then define by enlargement the sum operation $\sum $ for a
hyperfinite set of indices. This internal operation has by transfer the
following characteristic properties:

\begin{enumerate}
\item $\sum_{b\in \varnothing }\lambda _{b}b=0$

\item $\sum_{b\in B_{1}\cup B_{2}}\lambda _{b}b+\sum_{b\in B_{1}\cap
B_{2}}\lambda _{b}b=\sum_{b\in B_{1}}\lambda _{b}b+\sum_{b\in B_{2}}\lambda
_{b}b$, for $B_{1}$ and $B_{2}$ hyperfinite sets.
\end{enumerate}

\noindent This leads to:

\begin{proposition}
\label{prop5}Let $B$ be an orthonormal hyperbasis of a hyperhermitian space $%
F$. Every $x\in F$ could be expressed as: $x=\sum_{b\in B}<b,x>b$
\end{proposition}

\begin{remark}
\label{rem6}Infinite hyperfinite sets are uncountable.
\end{remark}

\subsection{Hull of a hyperhermitian space}

Let $F$ be a hyperhermitian space with a scalar product $<,>$. We define a
hypernorm $\left\Vert {}\right\Vert $(i.e. a norm with hyperreal values) by $%
\left\Vert x\right\Vert =\sqrt{<x,x>}$.

Let us consider now the set $F_{b}\subset F$ of bounded hypernorm elements
and define over $F_{b}$ an equivalence relation $\thicksim $ by : $%
x\thicksim y$ iff $x-y$ has an infinitesimal norm. The operations over $F$
are preserved by passing to the quotient space and make of $%
^{e}F=F/\thicksim $ a complex vector space. The scalar product $<,>$ is also
preserved by passing to the quotient and make of $^{e}F$ a prehilbertian
space. A classic result (see \cite{LuxStroy}) shows that $^{e}F$ is a
Hilbert space. For $x\in $ $^{e}F$, we call shadow of $x$ the class of $x$
in $^{e}F$\ and denote it also $^{\circ }x$.

Let $B$ be an orthonormal hyperbasis of $F$. If $F$ is a non
finite-dimensional space, set $B$ is uncountable, while for every bounded $x$%
, the subset $\{b\in B$ $/$ $<b,x>$ is not negligible$\}$ is countable. This
prevents the extension of the equality $x=\sum_{b\in B}<x,b>$\ to the hull
of $F$.

\subsection{The Loeb spectral measure\label{sec2.4}}

Let $F$ be an infinite-dimensional hyperhermitian space and $B$ an
orthonormal basis. The set $\mathfrak{T}_{0}(B)$ of hyperfinite subsets of $B
$ is a Boolean algebra where the union of a increasing sequence $(B_{n})$
does not belong to the Boolean algebra unless this sequence is stationary.
For every $\xi \in F_{b}$ we define a positive measure on $\mathfrak{T}%
_{0}(B)$ by $\nu _{\xi }(C)=$ $^{\circ }\left( \sum_{b\in C}\left\vert
<b,\xi >\right\vert ^{2}\right) $. By the fundamental theorem of P. Loeb,
this measure extends to the Boolean $\sigma $-algebra $\mathfrak{T}(B)$\
generated by $\mathfrak{T}_{0}(B)$ and depends only on the shadow $h$ of $%
\xi $.

\begin{theorem}
\label{thm7}There is a unique mapping $\varpi $ from the Boolean $\sigma $%
-algebra $\mathfrak{T}(B)$ to the set of orthogonal projectors of $^{e}F$
such that for any $h\in $ $^{e}F$ and any $C\in \mathfrak{T}(B):$ $%
v_{h}(C)=<h,\varpi (C)h>$.
\end{theorem}

\subsection{Hull of a symmetric internal linear function over a
hyperhermitian space}

Let $v$ be this function. The linear mapping $v-iI$ (where $I$ is the
identity function) is bijective and allows us to define the Cayley transform 
$c_{v}=(v+iI)(v-iI)^{-1}$\ which is a unitary function that becomes, by
passing to the quotient, a unitary function $^{e}c_{v}$\ over $^{e}F$. Let $%
^{e}F_{v,f}$ be the orthogonal of the eigenspace of $^{e}c_{v}$ associated
to the eigenvalue $1$.

The graph $\Gamma _{v}\subset F\times F$ is an internal subset. The shadows
of $\Gamma _{v}$ bounded elements in $^{e}(F\times F)$ make up a closed
subspace $G_{v}$. Note here that $^{e}(F\times F)$\ could be easily
identified with $^{e}F\times ^{e}F$.

\begin{proposition}
\label{prop8}The intersection $G_{v}\cap $ $^{e}F_{v,f}\times $ $^{e}F_{v,f}$%
\ is the graph of a selfadjoint operator $t_{v}$ in $^{e}F_{v,f}$.
\end{proposition}

Let us define an orthonormal hyperbasis $B$ of $F$ made of $v$ eigenvectors.
The mapping $b\mapsto \lambda _{b}$\ (which associates to an eigenvector $b$
its eigenvalue) is an internal hyperreal function. The set $B_{v,f}$ of $v$
eigenvectors with bounded eigenvalues is a borelian set in the sense of \ref%
{sec2.4}

\begin{proposition}
\label{prop9}The orthogonal projector over $^{e}F_{v,f}$ is the orthogonal
projector associated to the borelian set $B_{v,f}$ by the spectral measure $%
\varpi $ defined in \ref{sec2.4}
\end{proposition}

For any $\lambda \in \mathbb{R}$ the set $\Omega _{\lambda }=\left\{ b\in
B_{v,f}\mbox{ }|\mbox{ }^{\circ }\lambda _{b}\leq \lambda \right\} $, is a
borelian set. We can thus define $E_{\lambda }=\varpi (\Omega _{\lambda })$
as an orthogonal projector in $^{e}F_{v,f}$.

\begin{proposition}
\label{prop10}We have $t_{v}=\int_{-\infty }^{+\infty }\lambda dE_{\lambda }$%
, which means that $(E_{\lambda })$ is the spectral family of $t_{v}$.
\end{proposition}

\section{Renormalization of an orthonormal basis}

We still consider a hyperhermitian space $F$ and an orthonormal hyperbasis $%
B $ of $F$.

\subsection{$^{e}F_{B,h}$ Spaces}

For any element $h\in ^{e}F$, $h=^{\circ }x$ and $^{e}F_{B,h}$ denotes the
closed subspace generated by $\varpi (C)(h),$ where $C$ evolves in the Borel
algebra $\mathfrak{T}(B)$.

$\nu _{h}$ is the unique measure over $\mathfrak{T}(B)$ such that 
\begin{equation*}
\nu _{h}(C)=\text{ }^{\circ }\left( \sum_{b\in C}\left\vert <b,x>\right\vert
^{2}\right)
\end{equation*}%
for every hyperfinite subset $C$.

\begin{proposition}
\label{prop11}

\begin{description}
\item {a.} There is a unique isometry $J_{B,h}:\,^{e}F_{B,h}\rightarrow
L^{2}(B,\nu _{h})$ that links $\varpi (C)(h)$ to the characteristic function
class of $C$ in $L^{2}(B,\nu _{h})$.

\item {b.} Let $z$ be a bounded element of $F$ and $k$ its shadow in $%
^{e}F_{B,h}$. The image $\overline{k}$ of $k$ is the class of the function $%
f_{z}:b\mapsto ^{\circ }\left( \frac{<b,z>}{<b,x>}\right) $ if this ratio is
bounded and $0$ otherwise.

\item {c.} Let $\zeta $ be the measure $\nu _{k}-\left\vert \overline{k}%
\right\vert ^{2}\nu _{h}$. $\zeta $ and $\nu _{h}$\ are mutually singular
(and $\nu _{k}=\left\vert \overline{k}\right\vert ^{2}\nu _{h}+\zeta $ is
the Lebesgue decomposition of $\nu _{k}$\ with respect to $\nu _{h}$).
\end{description}
\end{proposition}

\begin{remark}
\label{rem12}

\begin{description}
\item {a.} The quantity $\frac{<b,z>}{<b,x>}$ is in general the ratio of two
infinitesimals. We cannot indeed express it through the shadows $h$ and $k$
of $x$ and $z$.

\item {b.} We will determine $J_{B,h}^{-1}$ in the next paragraph (Theorem %
\ref{thm14}c).
\end{description}
\end{remark}

\subsection{Introduction of a second scalar product\label{sec3.2}}

Let us consider a second internal scalar product $\left[ ,\right] $ over the
hyperhermitian space $F$. There exists a symmetric positive definite
internal linear function $j$ such that $\left[ x,y\right] =<j(x),j(y)>$. The
dual scalar product $\left\{ ,\right\} $ is defined by $\left\{ x,y\right\}
=<j^{-1}(x),j^{-1}(y)>$ and we have $\left\vert <x,y>\right\vert ^{2}\leq %
\left[ x,x\right] \left\{ y,y\right\} $ for every $x,y$.

We now introduce the mapping :

\begin{equation*}
\begin{tabular}{lll}
$^{e}F_{\left[ ,\right] }\times ^{e}F_{\left\{ ,\right\} }$ & $\rightarrow $
& $\mathbb{C}$ \\ 
$\left( a,b\right) $ & $\mapsto $ & $\left\langle a,b\right\rangle $%
\end{tabular}%
\end{equation*}
where : $\left\langle ^{\circ }x_{\left[ ,\right] }|^{\circ }y_{\left\{
,\right\} }\right\rangle =$ $^{\circ }<x,y>,$\newline
$^{e}F_{\left[ ,\right] }$ is the hull of $F$ with a scalar product $[,],$%
\newline
$^{e}F_{\left\{ ,\right\} }$ is the hull of $F$ with a scalar product $%
\{,\}, $\newline
$^{\circ }x_{\left[ ,\right] }$ is the shadow of $x$ in $^{e}F_{\left[ ,%
\right] },$\newline
$^{\circ }y_{\left\{ ,\right\} }$ is the shadow of $y$ in $^{e}F_{\left\{
,\right\} }$.\newline

\noindent\textbf{Reminder} - Let $G$ be a complex vector space. A mapping $%
f:G\rightarrow \mathbb{C}$ is called antilinear if $f(x+y)=f(x)+f(y)$ and $%
f(\lambda x)=\overline{\lambda }f(x)$. This property is equivalent to have $%
\overline{f}:x\mapsto \overline{f(x)}$ linear.

A linear function $f$ over a normed space $G$ is bounded, if the linear form 
$\overline{f}$\ is bounded. The norm of $f$ is then defined as the norm of $%
\overline{f}$. The vector space of bounded antilinear functions from $G$ to $%
\mathbb{C}$, is called the antilinear dual of $G$ and denoted $\widetilde{G}$%
.

\begin{proposition}
\label{prop13}Every linear form $f$ over $^{e}F_{\left\{ ,\right\} }$ can be
defined as \newline
$y\mapsto \left\langle x_{f},y\right\rangle $\ with a unique $x_{f}\in $ $%
^{e}F_{\left[ ,\right] }$. The mapping $f\mapsto x_{f}$\ is an antilinear
bijective isometry from the dual of $^{e}F_{\left\{ ,\right\} }$ into $%
^{e}F_{\left[ ,\right] }$.

Every antilinear form $\phi $ over $^{e}F_{\left[ ,\right] }$ can be defined
as $x\mapsto \left\langle x|y_{\phi }\right\rangle $ with a unique $y_{\phi
}\in $ $^{e}F_{\left\{ ,\right\} }$. The mapping $y_{\phi }\mapsto \phi $ is
an antilinear bijective isometry from the dual of $^{e}F_{\left[ ,\right] }$
into $^{e}F_{\left\{ ,\right\} }$.

In other words, the sesquilinear form $\left\langle |\right\rangle $ puts in
duality $^{e}F_{\left[ ,\right] }$ and $^{e}F_{\left\{ ,\right\} }$.
\end{proposition}

In addition we can find an orthonormal hyperbasis $\Gamma $\ of $F,[,]$\
which is also orthogonal for $<,>$.

The hypersum $\sigma =\sum_{\gamma \in \Gamma }\left\Vert <\gamma ,\gamma
>\right\Vert =\sum_{\gamma \in \Gamma }\sqrt{<\gamma ,\gamma >}$\ is
independent of the choice of $\Gamma $ basis .

\begin{theorem}
\label{thm14}Assume $\sigma $ bounded.

\begin{description}
\item {a.} The identity mapping from $F$ to $F$\ transforms bounded elements
for $[,]$ into bounded elements for $<,>$ and bounded elements for $<,>$
into bounded elements for $\{,\}$.This induces, by passing to the quotient,
continuous mappings : $^{e}F_{\left[ ,\right] }\rightarrow $ $^{e}F$ and $%
^{e}F\rightarrow $ $^{e}F_{\left\{ ,\right\} }$.

\item {b.} The complementary of the set $\left\{ b|\frac{b}{<x,b>}%
\mbox{ is
bounded for }\left\{ ,\right\} \right\} $ is of null $\nu _{k}$-measure. Let
us define the following function :%
\begin{equation*}
\begin{array}{cccc}
\phi _{x}: & B \rightarrow & ^{e}F_{\left\{ ,\right\} } &  \\ 
& b \mapsto &  & \left\{ 
\begin{array}{cc}
^{\circ }\left( \frac{b}{<x,b>}\right) & \mbox{ if }\frac{b}{<x,b>}%
\mbox{ is
bounded} \\ 
0 & \mbox{ if not}%
\end{array}%
\right.%
\end{array}%
\end{equation*}

\item {c.} Let $f\in L^{2}(B,\nu _{h})$ and $J_{B,h}:$ $^{e}F_{B,h}%
\rightarrow L^{2}(B,\nu _{h})$ the isometry defined in proposition \ref%
{prop11}. The mapping :%
\begin{equation*}
\begin{array}{cc}
B\rightarrow & ^{e}F_{\left\{ ,\right\} } \\ 
b\mapsto & f(b)\phi _{x}(b)%
\end{array}%
\end{equation*}%
is weakly integrable. For every $q\in $ $^{e}F_{h}$ we have : 
\begin{equation*}
<\widetilde{q},J_{B,h}^{-1}(f)>=\int_{B}f(b)\left\langle q,\phi
_{x}(b)\right\rangle d\nu _{h}(b)
\end{equation*}%
where $\widetilde{q}$ is the image of $q$ by the mapping : $^{e}F_{\left[ ,%
\right] }\rightarrow $ $^{e}F$ defined in a.
\end{description}
\end{theorem}

\begin{remark}
\label{rem15}

\begin{description}
\item a. Point c. could be also expressed as :\ the weak integral $%
\int_{B}f(b)\left\langle q,\phi _{x}(b)\right\rangle d\nu _{h}(b)$ is the
image of $J_{B,h}^{-1}(f)$\ by the mapping : $^{e}F\rightarrow
^{e}F_{\left\{ ,\right\} }$\ defined in point a.

\item b. With a possible exception for a countable set of elements $b$, the
scalar product $<x,b>$ is infinitesimal and $\frac{b}{<x,b>}$ is thus
unbounded for \newline
$<,>$. The introduction of the $\{,\}$ associated norm (renormalization
procedure) allows however to extract a finite value.
\end{description}
\end{remark}

\section{Application to spectral theory}

Let $G$ be a complex vector space with a locally convex topology, $<,>$ a
continuous scalar product over $G$, and $u$ a continuous linear function
over $G,$ symmetric and essentially selfadjoint relative to $<,>$. We have
the following results :

\begin{description}
\item {a.} $<u(x),y>=<x,u(y)>$ for every $x$ and $y$ in $G$,

\item {b.} If $\widehat{G}$\ is the completion of $G$ for the scalar product 
$<,>$, the closure of the graph of $u$ in $\widehat{G}\times \widehat{G}$ is
the graph of a selfadjoint operator $\overline{u}$.
\end{description}

\noindent \textbf{Note}

\begin{enumerate}
\item A scalar product $[,]$ over $G$ is nuclear (relative to $<,>$) if
there is a constant $C$ such that $\sum_{i=1}^{n}\left\Vert b_{i}\right\Vert
^{2}=\sum_{i=1}^{n}<b_{i},b_{i}>\leq C$\ for any finite sequence $%
b_{1,\cdots ,}b_{n}$ of vectors of $G$ that is orthonormal for $[,]$.

\item $[,]$ is called strongly nuclear (relative to $<,>$) if there is an
auxiliary nuclear scalar product $\left[ ,\right] _{aux}$ such that $[,]$ is
nuclear relative to $\left[ ,\right] _{aux}$.
\end{enumerate}

In what follows we assume a fixed scalar product $[,]$ over $G$ continuous
and strongly nuclear relative to $<,>$. The existence of such a scalar
product is guaranteed if $G$ is a nuclear space.

We introduce the enlargement $^{\ast }G$ of $G$ and provide it with the
internal scalar product obtained by extending the scalar product $<,>$ (with
the same notation). We then extend $u$ into an internal linear application $%
^{\ast }u$ over $^{\ast }G$. By proposition \ref{prop3} there exists a
hyperfinite-dimensional internal subspace $F$ of $^{\ast }G$ that contains
the image of $G$ in $^{\ast }G$ by the mapping $\ast $.\newline
Let us fix such a subspace $F$ that we assume hyperhermitian for the scalar
product $<,>$. To simplify notations, we identify $G$ with its image in $F$. 
\newline
The canonical injection of $G$ into $F$ induces an isometry from $G$ into $%
^{e}F$; and the closure $\overline{G}$ of $G$ in $^{e}F$ is a completion of $%
G$ for $<,>$.

Let $v$ be a mapping over $F$ that associates to an element $x$ the
orthogonal projection of $u(x)$ over $F$. $v$ is an internal linear
application, symmetric for $<,>$, that extends $u$. We can thus apply to $v$
propositions \ref{prop8} and \ref{prop9} using same notations :

\begin{proposition}
\label{prop16}

\begin{description}
\item {a.} The mapping $t_{v}$ associated to $v$ in proposition \ref{prop8}
is defined everywhere in $G$; and $t_{v}=u(x)$.

\item {b.} Let $(E_{\lambda })$ be the spectral family associated to $t_{v}$
in proposition \ref{prop10}. Every $E_{\lambda }$ preserves the closure $%
\overline{G}$ and the family $\left( E_{\lambda |\overline{G}}\right) $ of $%
E_{\lambda }$ restrictions to $\overline{G}$ is the spectral family of
projectors of the closure $\overline{u}$\ of $u$.
\end{description}
\end{proposition}

\begin{remark}
\label{rem17}Passing from $G,u$ to $F,v$ can be seen as the nonstandard
variant of the Ritz projection method.
\end{remark}

\subsection{Disintegration of the measure\ $\protect\nu _{h}$}

Let us fix an orthonomal hyperbasis $B$ of eigenvectors of $v$ and denote $%
\lambda _{b}$ the eigenvalue associated to eigenvector $b$ and $\varpi $ the
Loeb spectral measure associated to $B$ in theorem \ref{thm7}. Let $h$ be an
element of $\overline{G}$. Propositions \ref{prop9} and \ref{prop16} show
that $\sigma _{h}(\lambda )=<h,E_{\lambda }h>$ is equal to $<h,\varpi
(\Omega _{\lambda })(b)>$ where $\Omega _{\lambda }$ is the set of $b$ such
that $\lambda _{b}$ is bounded and $^{\circ }\lambda _{b}\leq \lambda $.

\begin{proposition}
\label{prop18}Let $q$ be a real function defined on $B$ as follows :%
\begin{equation*}
q(b)=\left\{ 
\begin{array}{cc}
^{\circ }\lambda _{b} & \mbox{if }\lambda _{b}\mbox{ is limited} \\ 
0 & \mbox{otherwise}%
\end{array}%
\right. 
\end{equation*}%
The Stieljes measure generated by the increasing real function $\sigma _{h}$%
\ is the image of measure $\nu _{h}$ by $q$.
\end{proposition}

Since a linear raising $j:L^{\infty }(\mathbb{R},d\sigma _{\lambda
})\rightarrow \mathcal{L}^{\infty }(\mathbb{R},d\sigma _{\lambda })$ (i.e. a
linear mapping such as $\phi $ is the class of $j(\phi )$) exists, with $%
\mathcal{L}^{\infty }(\mathbb{R},d\sigma _{\lambda })$ the space of
essentially bounded numerical functions, we have the following proposition :

\begin{proposition}
\label{prop19}We can associate to every hyperfinite subset $C$ of $B$ an
element $\omega _{C}$ of the space of essentially bounded numerical
functions $\mathcal{L}^{\infty }(\mathbb{R},d\sigma _{\lambda })$, that
meets the following conditions :

\begin{description}
\item {a.} The{\ function} $\omega _{C}$\ takes only positive real values

\item {b.} $\omega _{\varnothing }(\lambda )=0$ and $\omega _{C_{1}\cup
C_{2}}(\lambda )=\omega _{C_{1}}(\lambda )+\omega _{C_{2}}(\lambda )$\ if $%
C_{1}$ and $C_{2}$ are disjoint.

\item {c.} Every numerical function $f$ integrable with respect to $d\sigma
_{\lambda }$ verifies: $\int_{C}f(q(b))d\nu _{h}(b)=\int_{-\infty }^{+\infty
}f(\lambda )\omega _{C}(\lambda )d\sigma _{h}(\lambda ).$
\end{description}
\end{proposition}

Statements a. and b. show that we can extend the mapping \newline
$\tau _{\lambda }:C\mapsto \omega _{C}(\lambda )$ to a Loeb measure over the
Boolean $\sigma $-algebra $\mathfrak{T}(B)$.

\begin{proposition}
\label{prop20}

\begin{description}
\item {a.} $\nu _{h}(C)=\int_{-\infty }^{+\infty }\tau _{\lambda }(C)d\sigma
_{h}(\lambda )$ for every borelian set $C$.

\item {b.} For $\sigma_{h}$-almost every $\lambda$, the measure $\tau
_{\lambda }$ is carried by the set\newline
$q^{-1}(\lambda )=\left\{ b\in B|\lambda_{b}\mbox{ is bounded and }^{\circ
}\lambda _{b}=\lambda \right\} $
\end{description}
\end{proposition}

This means that the family $(\tau _{\lambda })$\ is a disintegration of the
measure\ $\nu _{h}$ with respect to its image by $q$.

We assumed so far the existence of a scalar product $[,]$ that is strongly
nuclear over $G$. This scalar product extends to $^{\ast }G$ and induces a
scalar product on $F$ with hypercomplex values. We will keep denoting it $%
[,] $. We can now introduce the dual scalar product $\{,\}$ (see \ref{sec3.2}%
) :

\begin{lemma}
\label{lem21}

\begin{description}
\item {a.} Let $\Gamma $ be an orthonormal basis for $[,]$ that is also
orthogonal for $<,>$. The hypersum $\sum_{\gamma \in \Gamma }\left\Vert
\gamma \right\Vert $ is bounded -- and we can thus apply theorem \ref{thm14}.

\item {b.} Let $\phi $ be an element of $F$ limited for the dual scalar
product $\{,\}$. For every $g\in G$ the hypercomplex $<g,\phi >$ is limited
and the function $\phi ^{G}:g\mapsto ^{\circ }\left( <g,\phi >\right) $ is a
continuous antilinear form on $G$.
\end{description}
\end{lemma}

\noindent\textbf{Consequences of lemma \ref{lem21}}

\begin{description}
\item {a.} $\eta _{b}=\frac{b}{\sqrt{\left\{ b,b\right\} }}$ is of norm $1$
for the dual scalar product $\{,\}$ and defines therefore the continuous
antilinear form $\eta _{b}^{G}$\ on $G$.

\item {b.} Assume again theorem \ref{thm14} notations. If $h$ is an element
of $^{e}F$ and a shadow of an element $x$, we know that $\nu _{h}$-almost
every $\frac{b}{<x,b>}$ is limited for $\{,\}$ and we have therefore a
continuous antilinear form $\phi _{x}^{G}:g\mapsto ^{\circ }\left( \frac{%
<g,b>}{<x,b>}\right) $.
\end{description}

\noindent Consider an element $h$ of $G$. $h$ is its own shadow. Since 
\newline
$\frac{b}{<h,b>}=\frac{\frac{b}{\sqrt{\left\{ b,b\right\} }}}{<h,\frac{b}{%
\sqrt{\left\{ b,b\right\} }}>}$, $\frac{b}{<h,b>}$ is limited for $\{,\}$
iff $<h,\frac{b}{\sqrt{\left\{ b,b\right\} }}>$ is not infinitesimal (i.e.
if $\eta _{b}^{G}(h)\neq 0)$. And we have : $\phi _{b}^{G}(b)(g)=\frac{\eta
_{b}^{G}(g)}{\eta _{b}^{G}(h)}$.

\subsection{ Main Theorem}

Let $G_{a}^{\prime }$ be the antilinear dual of $G$ (i.e. the space of
continuous antilinear forms on $G$). $G_{a}^{\prime }$ is obviously a
complex subspace of the space of complex functions defined on $G$.\newline
The function $u$ admits an adjoint $u^{\ast }$ which is the linear function
from $G_{a}^{\prime }$ into itself defined by $u^{\ast }(\phi
)(g)=\phi(u(g)) $. We will call eigenfunctional of $u$ any element of $%
G_{a}^{\prime }$ which is an eigenvector of $u^{\ast }$. Those of $\eta
_{b}^{G}$ which are not null are such eigenfuctionals of $u$ (for the
eigenvalues $^{\circ}\lambda _{b}$).\newline

We finally get to the last and main theorem of this work (with the same
notations used so far) :

\begin{theorem}
\label{thm22}Let $h$ be an element of $G$ and $B_{h}$ the set of $b\in B$
such that $\eta _{b}^{G}(h)\neq 0$.

\begin{description}
\item {a.} The set $B_{h}$ is a borelian subset of $B$ whose complementary
is of null $\sigma _{h}$-measure. For $b\in B_{h}$ the antilinear form $\eta
_{b}^{G}$ is an eigenfunctional of $u$ for the eigenvalue $^{\circ }\lambda
_{b}$.

\item {b.} The family $\left( \phi _{h}^{G}(b)=\frac{\eta _{b}^{G}}{\eta
_{b}^{G}(h)}\right) ,b\in B_{h}$, is weakly $\tau _{\lambda }$-integrable
for $\sigma _{h}$-almost every $\lambda $.

\item {c.} For almost every $\lambda $\ the weak integral $\omega
_{h,\lambda }=\int_{B_{h}}\frac{1}{\eta _{b}^{G}(h)}\eta _{b}^{G}d\tau
_{\lambda }(b)$ is, if not null, an eigenfunctional of $u$ for the value $%
\lambda $.

\item {d.} Let $(E_{\lambda })$ be the spectral projector family of the
closure of $u$ in the completion of $G$ for $<,>$ (see proposition \ref%
{prop16}). For any $g\in G$ we have $<g,E_{\lambda }h>=\int_{-\infty
}^{\lambda }\omega _{\mu }(g)d\sigma _{h}(\mu )$
\end{description}
\end{theorem}

\bigskip 

\noindent \textbf{Variant} -- The antilinear form $g\mapsto <g,E_{\lambda }h>
$\ associated to an element $E_{\lambda }h$\ of the completion of $G$, is
the weak integral $\int_{-\infty }^{\lambda }\omega _{\mu }d\sigma _{h}(\mu )
$ where for almost every $\mu $ the antilinear form $\omega _{\mu }$ is, if
not null, an eigenfunctional of $u$ for the eigenvalue $\mu $ of $u$.

\begin{remark}
The statement c. of the preceding theorem is a nonstandard approach of the
integral representation of the spectral projectors of an essentially
selfadjoint operator with eigenfunctionals (see \cite{GelChil} Tome 3
Chapter 4 n${{}^{\circ }}$ 4.2 Theorem 1).\newline
The outcome of this approach is the description of these eigenfunctional $%
\omega _{t}$ with the eigenvectors of an extension $v$ of the considered
operator $u$ in hyperfinite dimension. These eigenvectors do not have an
immediate interpretation in the initial space $G$; but their normalised
forms for the scalar product $\{,\}$ define eigenfunctionals.\newline
We aggregate through integration those eigenfunctionals which are associated
to the same eigenvalue; and here again the Loeb's technique comes to make
things much easier.
\end{remark}

\noindent \textbf{Acknowledgements }

\noindent I would like to thank Professor Fran\c{c}ois Aribaud for his
helpful comments and suggestions. All errors and omissions remain mine.


\begin{thebibliography}{Luxemburg Stroyan}
\bibitem[Gelfand Chilov]{GelChil} - I.M Gelfand and G.E Chilov - Les
distibutions Tome 3 (Equations diff\'{e}rentielles) Dunod Paris 1965

\bibitem[Goldblatt]{Gold} - Robert Goldblatt - Lectures on the Hyperreals -
An Introduction to Nonstandard Analysis Graduate Texts in Mathematics n${%
{}^\circ}$ 188 - Springer Verlag 1998

\bibitem[Jech]{Jech} - Thomas Jech - Set Theory - Pure and Applied
Mathematics - Academic Press 1978

\bibitem[Loeb Wolff]{LW} - Peter A. Loeb and Manfred Wolff - Nonstandard
Analysis for the working mathematician Mathematics and its applications vol
510 Kluwer Academic Publishers 2000

\bibitem[Luxemburg Stroyan]{LuxStroy} - W.A.J Luxemburg and K.D.Stroyan -
Introduction to the Theory of Infinitesimals - Academic Press 1976
\end{thebibliography}
\end{document}